\documentstyle[amsfonts,11pt]{article}

\begin{document}
\title{Answering a question on relative countable paracompactness}
\author{M. V. Matveev\\
{\small Department of Mathematics, University of California, Davis,}\\
{\small Davis, CA, 95616, USA (address valid till June 30, 2000)}\\
{\small e-mail misha$\underline{\mbox{ }}$matveev@hotmail.com}
}
\date{}
\maketitle

In \cite{Yasui}, Yoshikazu Yasui formulates some results
on relative countable paracompactness and poses some questions.
Like it is the case with many other topological properties
\cite{Arh}, countable paracompactness has several
possible relativizations.
Thus a subspace $Y\subset X$ is called {\em countably 1-paracompact}
in $X$ provided for every countable open cover $\cal U$ of $X$
there is an open cover $\cal V$ of $X$ which refines $\cal U$ and is
locally finite at the points of $Y$ (i.e. every point of $Y$ has
a neighbourhood in $X$ that meets only finitely many
elements of $\cal V$).
Yasui asserts that if a countably compact space $Y$ is closed
in a normal space $X$ then $Y$ is countably 1-paracompact in $X$
and asks (Problem 1 in\cite{Yasui}) if normality can be omitted.
The answer is negative as it is demonstrated by
$X=((\omega_1+1)\times(\omega+1))\setminus\{(\omega_1,\omega)\}$, 
$Y=\omega_1\times\{\omega\}$ and
${\cal U}=\{\omega_1\times(\omega+1)\}\cup\{\{(\omega_1+1)
\times\{n\}\}:n\in\omega\}$.
This well-known construction provides also the following
general statement (recall that a space is Linearly Lindel\"of
iff every uncountable set of regular cardinality
has a complete accumulation point see e.g. \cite{AB}).

\bigskip
\noindent{\bf Theorem 1} {\em If a Tychonoff space $Y$
is countably 1-paracompact in every Tychonoff space $X$
that contains $Y$ as a closed subspace then $Y$ is 
linearly Lindel\"of.
}

\bigskip
\noindent{\bf Proof:}
Suppose not. 
Then there is an uncountable set $Z\subset Y$ of regular cardinality
and without complete accumulation points in $Y$.
Enumerate $Z$ as $\{z_\alpha:\alpha<\kappa\}$
where $\kappa=|Z|$
and put $Z^*=Z\cup\{z^*\}$ where $z^*\not\in Z$.
Further, put
$X=(Z^*\times\omega)\cup(Y\times\{\omega\})$.
Topologize $X$ as follows.
The points of $Z\times\omega$ are isolated.
A basic neighbourhood of a point $(z^*,n)$, where $n\in\omega$,
takes the form $\{(z_\gamma,n):\gamma>\alpha\}\cup\{(z^*,n)\}$
where $\alpha<\kappa$.
A basic neighbourhood of a point $(y,\omega)\in Y\times\{\omega\}$
takes the form 
$O_{Un}=((U\cap Z)\times\{m\in\omega:m>n\})\cup(U\times\{\omega\})$
where $U$ is a neighbourhood of $y$ in $Y$
and $n\in\omega$.
Then $X$ contains
a closed subspace $\tilde{Y}=Y\times\{\omega\}$
homeomorphic to $Y$.

Now we check that $X$ is a Tychonoff space.
It is clear that the points of $Z^*\times\omega$ have 
local bases consisting of clopen sets.
So let $x=(y,\omega)\in\tilde{Y}$
and let $O$ be a neighbourhood of $x$ in $X$.
Then there are a neighbourhood $U$ of $y$ in $Y$ and $n\in\omega$
such that $x\in O_{Un}\subset O$.
Further, there is a neighbourhood $V$ of $y$ in $Y$
such that $y\in V\subset U$ and $|V\cap Z|<\kappa$.
Since $Y$ is Tychonoff, there is a function $f:Y\to\Bbb R$
such that $f(y)=0$ and $f(Y\setminus V)=\{1\}$.
Define a function $\tilde{f}:X\to\Bbb R$ as
$$
\tilde{f}(u)=\left\{\begin{array}{ll}
f(v), & \mbox{ if }u=(v,\omega)\\
f(z), & \mbox{ if }u=(z,m)\mbox{ and }m>n\\
1, & \mbox{ otherwise.}
\end{array}\right.
$$
Then $\tilde{f}$ is continuous; this follows from the inclusion
$\overline{\tilde{f}^{-1}({\Bbb R}\setminus\{1\})}
\subset(Z\times\omega)\cup(Y\times\{\omega\})$.
Finally, $\tilde{f}(x)=0$ and $\tilde{f}(X\setminus U)=\{1\}$.
So $X$ is Tychonoff.

It follows from the inequality ${\rm cf}(\kappa)>\omega$ that
the open cover
${\cal U}=\{(Z\times\omega)\cup(Y\times\{\omega\})\}\cup
\{Z^*\times\{n\}:n\in\omega\}$ of $X$ does not
have an open, locally finite at all points 
of $\tilde{Y}$ refinement.
So $\tilde{Y}\sim Y$ is not countably 1-paracompact in $X$.
$\Box$

\bigskip
\noindent{\bf Theorem 2} {\em A Tychonoff countably compact
space $Y$ is countably 1-paracompact in every Tychonoff $X\supset Y$
iff $Y$ is compact.}

\bigskip
\noindent{\bf Proof:}
Necessity follows from the previous theorem and the fact that
every countably compact (in fact, even every countably paracompact,
see \cite{Mis}) linearly Lindel\"of space is compact.

Routinous proof of sufficiency is omitted.
$\Box$

\bigskip
\noindent{\bf Theorem 3}
{\em If a Lindel\"of space $Y$ is a closed subspace of a regular
space $X$ then $Y$ is countably 1-paracompact in $X$.}

\bigskip
\noindent{\bf Proof:}
Let ${\cal U}=\{U_n:n\in\omega\}$ be a countable open cover of $X$.
For every $y\in Y$ fix $n(y)\in\omega$ and an open set
$W_y\subset X$ so that
$y\in W_y\subset\overline{W_y}\subset U_{n(y)}$.
The cover ${\cal W}=\{W_y:y\in Y\}$
contains a countable subcover of $Y$, say
$\{W_{y_k}:k\in\omega\}$.
Then ${\cal V}=\{U_{n(y_k)}\setminus
\overline{\cup\{W_{y_l}:l<k\}}:k\in\omega\}\cup
\{U_n\setminus\overline{\cup\{W_{yl}:l<n\}}:n\in\omega\}$
is an open refinement of $\cal U$
and $\cal V$ is locally finite at all points of $Y$.
$\Box$

\bigskip
\noindent{\bf Remark 1}
In Theorem 1 and Theorem 2 one can replace ``Tychonoff''
with ``regular''

\bigskip
\noindent{\bf Remark 2}
The results above are similar to some results about normailty
and property (a) from \cite{BY}, \cite{MPT}.

\bigskip
The paper was written when the author was visiting University of California,
Davis. The author expresses his gratitude to colleagues
from UC Davis for their kind hospitality.


\end{document}